\documentclass[reqno]{amsart}

\usepackage{amssymb}
\usepackage[english]{babel}

\newtheorem{theorem}{Theorem}
\newtheorem{lemma}{Lemma}
\newtheorem{corollary}{Corollary}
\theoremstyle{definition}

\theoremstyle{remark}
\newtheorem{remark}{Remark}

\newcommand{\norm}[1]{\left\Vert#1\right\Vert}

\newcommand{\spaceR}{{\mathbb{R}}}
\newcommand{\spaceL}{{\mathbf L}}
\newcommand{\spaceLI}{{\mathbf L}_{\infty}}
\newcommand{\spaceC}{{\mathbf C}}
\newcommand{\spaceAC}{{\mathbf{AC}}}

\DeclareMathOperator*{\vraisup}{vrai\, sup}

\newcommand{\unit}{{\mathbf1}{\,}}
\newcommand{\opCL}{\spaceC\to\spaceL}
\newcommand{\opCLI}{\spaceC\to\spaceL_{\infty}}
\newcommand{\cal}{\mathcal}
\newcommand{\calT}{{\cal T}}
\newcommand{\calA}{{\cal A}}
\newcommand{\calB}{{\cal B}}
\renewcommand{\a}{{\cal T}^+}
\renewcommand{\b}{{\cal T}^-}
\newcommand{\w}{\widetilde}
\newcommand{\tri}{\triangle}

\begin{document}

\title[Conditions for the solvability of the Cauchy problem]
{Conditions for the solvability of the Cauchy problem for linear first-order functional differential equations}

\author{E.I. Bravyi}

\address{State National Research Po\-ly\-tech\-nic University of Perm, Perm, Russia}

\email{bravyi@perm.ru}



\date{17.06.2013}

\begin{abstract}
Conditions for the unique solvability of the Cauchy problem for a family of scalar functional differential equations are obtained. These conditions are sufficient for the solvability of the Cauchy problem for every equation from the family and are necessary for the solvability of the Cauchy problem for all equations from the family.  In contrast to many known articles, we consider equations with functional operators acting into the space of essentially bounded functions.
\end{abstract}

\maketitle

\section*{Introduction}

In this article, the well-known integral conditions for the solvability of the Cauchy problem for linear functional differential equations (Theorem \ref{t-1}) are added to
necessary and sufficient conditions with point-wise restrictions on functional operators (Theorem \ref{t-2}). Also some conditions for solvability of the Cauchy problem for a family of quasilinear equations are obtained.

We use the following notation:
$\spaceR$ is the space of real numbers;
$\spaceC=\spaceC[a,b]$ is the Banach space of continuous functions $x:[a,b]\to\spaceR$ with the norm
\begin{equation*}
\norm{x}_{\spaceC}=\max_{t\in[a,b]}|x(t)|;
\end{equation*}
$\spaceL_\infty=\spaceL_\infty[a,b]$ is the Banach space of essentially bounded measurable functions $z:[a,b]\to\spaceR$ with the norm
\begin{equation*}
\norm{x}_{\spaceLI}=\vraisup_{t\in[a,b]}|x(t)|;
\end{equation*}
$\spaceL=\spaceL[a,b]$
is the Banach space of integrable functions $z:[a,b]\to\spaceR$ with the norm
\begin{equation*}
\norm{z}_{\spaceL}=\int_a^b|z(t)|\,dt,
\end{equation*}
it is supposed that all inequalities and equalities with functions from $\spaceL$ and  $\spaceL_\infty$ hold almost everywhere on $[a,b]$;
$\spaceAC=\spaceAC[a,b]$ is the Banach space of absolutely continuous functions
$x:[a,b]\to\spaceR$ with the norm
\begin{equation*}
\norm{x}_{\spaceAC}=|x(a)|+\int_a^b|\dot x(t)|\,dt,
\end{equation*}
$\unit(t)\equiv1$ is the unit function;
an operator $T:\opCL$ is said to be positive (or isotonic in the terminology of \cite{AMR})  if it maps each non-negative continuous function into an
almost everywhere non-negative function.

Consider the Cauchy problem for a first-order functional differential equation
\begin{gather}
\dot x(t)=(T^+x)(t)-(T^-x)(t)+f(t),\quad t\in[a,b], \label{e-1}\\
x(a)=c, \label{e-2}
\end{gather}
where $T^+$, $T^-:\opCL$ are linear positive operators, $f\in\spaceL$, $c\in\spaceR$.

A solution of \eqref{e-1}--\eqref{e-2} is a satisfying the initial conditions \eqref{e-2} function $x\in\spaceAC$ such that equality \eqref{e-1} holds almost everywhere on $[a,b]$.
Problem \eqref{e-1}--\eqref{e-2} is called uniquely solvable if it has a unique solution for every pair $f\in\spaceL$,  $c\in\spaceR$.

The positiveness of operators $T^+$ and $T^-$ implies their $u$-boundedness   (or the strong boundedness in other terminology) \cite{Zabrejko}. This property guarantees the Fredholm property of problem \eqref{e-1}--\eqref{e-2} \cite{AMR, Maksimov1}. From the Fredholm property it follows that problem \eqref{e-1}--\eqref{e-2} is uniquely solvable if and only if the problem has a unique solution for at least one pair  $f\in\spaceL$, $c\in\spaceR$. In particular, the problem is uniquely solvable if and only if the homogeneous problem
\begin{gather}\label{e-3}
\left\{
  \begin{array}{l}
\dot x(t)=(T^+x)(t)-(T^-x)(t),\quad t\in[a,b], \\
x(a)=0,
  \end{array}
\right.
\end{gather}
has only the trivial solution.

If linear positive operators $T^+$, $T^-:\opCL$ are Volterra type operators, then the Cauchy problem is uniquely solvable without additional conditions (see, for example, \cite{AMR}).
Some optimal, in a sense, conditions for the unique solvability of problem \eqref{e-1}--\eqref{e-2} are found in \cite{Hakl} (see also \cite{HL2002}, \cite{HLP2002}, \cite{HLS2002-5}) for, generally speaking, non-Volterra operators.
We give here this result in the form of necessary and sufficient conditions for the solvability. Note, that each linear positive operator $T:\opCL$ is bounded, its norm is defined by the equality
\begin{equation*}
    \norm{T}_{\opCL}=\int_a^b (T\unit)(t)\,dt.
\end{equation*}

\begin{theorem}[\cite{HL2002}]\label{t-1}
Let non-negative numbers $\calT^+$, $\calT^-$ be given. Problem \eqref{e-1}--\eqref{e-2} is uniquely solvable
for all linear positive operators $T^+$, $T^-:\opCL$ with the given norms $\norm{T^+}_{\opCL}=\calT^+$, $\norm{T^-}_{\opCL}=\calT^-$
if and only if
the inequalities
\begin{equation}\label{e-4}
\calT^+<1,\quad \calT^-<1+2\sqrt{1-\calT^+}
\end{equation}
hold.
\end{theorem}

As far as we know, unimprovable conditions for the solvability of the Cauchy problem \eqref{e-1}--\eqref{e-2}
in terms of norm operators $T^+$, $T^-:\spaceC\to\spaceLI$ are not obtained yet.

The main result is the statement, which is similar to Theorem \ref{t-1}, but deals with operators acting from the space $\spaceC$ into the space os essentially bounded functions $\spaceLI$. The norm of linear positive operator $T:\opCLI$ is defined by the equality
\begin{equation*}
    \norm{T}_{\opCLI}=\vraisup_{t\in[a,b]}\, (T\unit)(t).
\end{equation*}

For short, we use the notation
\begin{equation*}
    \calA\equiv (b-a)\,\a,\quad \calB\equiv (b-a)\,\b.
\end{equation*}

From Theorem \ref{t-1} it is easily to archive a sufficient condition for the solvability.
\begin{corollary}\label{c-1}
Let non-negative numbers $\calT^+$, $\calT^-$ be given. Then, for problem \eqref{e-1}--\eqref{e-2} to be uniquely solvable for all linear positive operators $T^+$, $T^-:\opCLI$ with given norms  $\norm{T^+}_{\opCLI}=\calT^+$, $\norm{T^-}_{\opCLI}=\calT^-$  it is necessary and sufficient that
\begin{equation}\label{e-5-new}
\calA<1,\quad    \calB< 1+2\sqrt{1-\calA}.
\end{equation}
\end{corollary}

The forthcoming Theorem \ref{t-2} shows that for all $\calA<1$ (except $\calA=0$)
conditions  \eqref{e-5-new} can be improved, and for $\calA=0$ the necessary and sufficient condition for the solvability of the problem
\begin{gather*}
\dot x(t)=-(T^-x)(t)+f(t),\quad t\in[a,b], \\
x(a)=c,
\end{gather*}
for all linear positive operators
$T^-:\opCLI$ with given norm $\norm{T^-}_{\opCLI}=\calT^-$
remains the inequality $\calB<3$ as in Theorem \ref{t-1}.

\section{The main results}

\begin{theorem}\label{t-2}
Let non-negative numbers $\calT^+$, $\calT^-$ be given. Then, for problem \eqref{e-1}--\eqref{e-2} to be uniquely solvable for all linear positive operators
$T^+$, $T^-:\opCLI$ with given norms $\norm{T^+}_{\opCLI}=\calT^+$, $\norm{T^-}_{\opCLI}=\calT^-$, it is necessary and sufficient that
\begin{equation}\label{e-13}
   (\calB^2-\calA^2)\,t^2+(\calA^2-\calB^2 +\calB)\,t+1-\calA>0
\end{equation}
for all $t\in[0,1]$.
\end{theorem}

\begin{corollary}
Let non-negative numbers $\calT^+$, $\calT^-$ be given. Then, for problem \eqref{e-1}--\eqref{e-2} to be uniquely solvable for all linear positive operators
$T^+$, $T^-:\opCLI$ with given norms $\norm{T^+}_{\opCLI}=\calT^+$, $\norm{T^-}_{\opCLI}=\calT^-$, it is necessary and sufficient that at least one of the following conditions be fulfilled:
\begin{enumerate}
\item[1)]
\begin{equation}\label{e-5}
\calA<1,
\end{equation}
\begin{gather}\label{e-6}
\begin{array}{l}
(2\calB-\calA)^2-\calA^2-
           \left(\calB^2-\calA^2-\calB+2\calA\right)^2>0\
\text{ for }\  \calB> \dfrac{1+\sqrt{1+\calA^2}}{2};
\end{array}
\end{gather}
\item[2)]
\begin{equation}\label{e-14}
\calA<1,\quad\calB<\min_{t\in(0,1)}\dfrac{t+\sqrt{(2t(1-t)\calA-1)^2+(1-t)(3t-1)}}{2t(1-t)};
\end{equation}
\item[3)]
\begin{equation*}
    \calA<1\ \text{ for } \calB\in\left[0,(1+\sqrt{5})/{2}\right],
\end{equation*}
\begin{equation*}
    \calA<\min_{t\in(0,1)}\dfrac{1-\sqrt{(2t(1-t)\calB-t)^2-(1-t)(3t-1)}}{2t(1-t)}\quad \ \text{ for }\ \calB\in\left((1+\sqrt{5})/{2},3\right).
\end{equation*}
\end{enumerate}
\end{corollary}

Let us remark an obvious corollary of Theorem \ref{t-2}, which can be easily proved by the Schauder fixed point theorem (see, for example, \cite[p.~190]{HLS2002-5}).
Consider the quasilinear Cauchy problem
\begin{gather}\label{e-33}
\left\{
  \begin{array}{l}
\dot x(t)=(T^+x)(t)-(T^-x)(t)+(Fx)(t),\quad t\in[a,b], \\
x(a)=c,
  \end{array}
\right.
\end{gather}
where $T^+$, $T^-:\opCLI$ are linear positive operators, the operator $F:\opCLI$ is continuous and bounded (maps bounded sets into bounded ones), $c\in\spaceR$.

\begin{corollary}
Suppose that linear positive operators $T^+$, $T^-:\opCLI$ satisfy the conditions of Theorem \ref{t-2}, a continuous bounded operator
$F:\opCLI$ satisfies the under linear growth condition
\begin{equation*}
    \lim\limits_{\norm{x}_{\spaceC}\to\infty} \frac{\norm{Fx}_{\spaceLI}}{\norm{x}_{\spaceC}}=0.
\end{equation*}
Then the Cauchy problem \eqref{e-33} has a solution.
\end{corollary}

\section{Proofs}
To prove Theorem \ref{t-2} and Corollary \ref{c-1} we need auxiliary assertions.

\begin{lemma}\label{l-1}
Let non-negative functions $p^+$, $p^-\in\spaceLI$ be given. To problem \eqref{e-1}--\eqref{e-2} to be uniquely solvable for all linear positive operators
$T^+$, $T^-:\opCLI$ satisfying the equalities $T^+\unit=p^+$ and $T^-\unit=p^-$,
it is necessary and sufficient that this problem be uniquely solvable for all linear positive operators
$T^+$, $T^-:\opCL$ such that  $T^+\unit\le p^+$, $T^-\unit\le p^-$.
\end{lemma}

\proof
It is clear that only the assertion on the necessity is needed to be proven.
Suppose that problem \eqref{e-1}--\eqref{e-2} is not uniquely solvable for some linear positive operators
$T^+$, $T^-:\opCLI$ such that $T^+\unit\le p^+$, $T^-\unit\le p^-$, therefore \eqref{e-3}
has a non-trivial solution. Then for perturbated operators
\begin{eqnarray*}
(\w T^+x)(t)&    \equiv & (T^+x)(t)+\left(p^+(t)-(T^+\unit)(t)\right)\,x(a),\quad t\in[a,b], \\
(\w T^-x)(t)&    \equiv & (T^-x)(t)+\left(p^-(t)-(T^-\unit)(t)\right)\,x(a),\quad t\in[a,b],
\end{eqnarray*}
which satisfy the conditions
\begin{equation*}
\w T^+\unit=p^+,\quad \w T^-\unit=p^-,
\end{equation*}
the homogeneous problem
\begin{gather*}
\left\{
  \begin{array}{l}
\dot x(t)=(\w T^+x)(t)-(\w T^-x)(t),\quad t\in[a,b], \\
x(a)=0,
  \end{array}
\right.
\end{gather*}
has the same non-trivial solution. Hence, the corresponding non-ho\-mo\-ge\-ne\-ous problems are not uniquely solvable.
\endproof

Taking in Lemma \ref{l-1} $p^+(t)=\calT^+$, $p^-(t)=\calT^-$, $t\in[a,b]$, we obtain the following result.

\begin{corollary}
Let non-negative numbers  $\calT^+$, $\calT^-$ be given. For problem \eqref{e-1}--\eqref{e-2} to be uniquely solvable for all linear positive operators
$T^+$, $T^-:\opCLI$ with given norms
\begin{equation*}
\norm{T^+}_{\opCLI}=\calT^+,\quad \norm{T^-}_{\opCLI}=\calT^-,
\end{equation*}
it is necessary and sufficient that the problem be uniquely solvable for all linear positive operators
$T^+$, $T^-:\opCLI$ such that
\begin{equation*}
(T^+\unit)(t)=\calT^+,\quad (T^-\unit)(t)=\calT^-,\quad t\in[a,b].
\end{equation*}
\end{corollary}

\begin{lemma}\label{l-2}
Let non-negative functions $p^+$, $p^-\in\spaceLI$ be given. For problem \eqref{e-1}--\eqref{e-2} to be uniquely solvable for all linear positive operators
 $T^+$, $T^-:\opCLI$ such that
\begin{equation}\label{e-7}
T^+\unit=p^+,\quad T^-\unit=p^-,
\end{equation}
it is necessary and sufficient that the problem
\begin{gather}\label{e-8}
\left\{
  \begin{array}{l}
\dot x(t)=p_1(t)x(\tau_1)+p_2(t)x(\tau_2),\quad t\in[a,b], \\
x(a)=0,
  \end{array}
\right.
\end{gather}
to have only the trivial solutions for all points $\tau_1$, $\tau_2\in[a,b]$ and for all functions $p_1,p_2\in\spaceLI$ satisfying the conditions
\begin{equation}\label{e-9}
    p_1+p_2=p^+-p^-,\quad -p^-\le p_1 \le p^+.
\end{equation}
\end{lemma}

\proof
To prove the sufficiency suppose that for given non-negative functions $p^+$, $p^-\in\spaceLI$ there exist linear positive operators
$T^+, T^-:\opCLI$ such that
\begin{equation*}
T^+\unit=p^+,\quad T^-\unit=p^-
\end{equation*}
and problem \eqref{e-1}--\eqref{e-2} is not uniquely solvable.
Then there exists a non-trivial solution $y\in\spaceAC$ of the homogeneous problem \eqref{e-3}.
Let $\tau_1$ be a point of the minimum, $\tau_2$ be a point of the maximum of the solution $y$. Then
\begin{gather*}
y(\tau_1)\,p^+(t)=    y(\tau_1)\,(T^+\unit)(t)\le (T^+ y)(t) \le     y(\tau_2)\,(T^+\unit)(t)=y(\tau_2)\,p^+(t),\ t\in[a,b],\\
y(\tau_1)\,p^-(t)=     y(\tau_1)\,(T^-\unit)(t)\le (T^- y)(t) \le     y(\tau_2)\,(T^-\unit)(t)=y(\tau_2)\,p^-(t),\ t\in[a,b].
\end{gather*}
Therefore,
\begin{equation*}
    y(\tau_1)\,p^+-y(\tau_2)\,p^-\le     T^+y-T^-y
\le     y(\tau_2)\,p^+-y(\tau_1)\,p^-.
\end{equation*}
It follows that there exists a measurable function $\xi:[a,b]\to[0,1]$ such that for the functions
\begin{equation*}
 p_1\equiv\xi\, p^+-(1-\xi)\,p^-,\quad     p_2\equiv(1-\xi)\, p^+-\xi\, p^-,
\end{equation*}
the equality
\begin{equation*}
(T^+y)(t)-(T^-y)(t)=p_1(t)\,y(\tau_1)+p_2(t)\,y(\tau_2),\quad t\in[a,b],
\end{equation*}
holds. It is clear that conditions \eqref{e-9} for the functions $p_1$, $p_2$ are fulfilled and problem \eqref{e-8} has a non-trivial solution.

Let us prove the necessity. Suppose that function $p^+$, $p^-\in\spaceLI$ are non-negative. Let conditions \eqref{e-9} be fulfilled and problem  \eqref{e-8} have a non-trivial solution.
Define linear positive solutions  $T^+$, $T^-:\opCLI$ by the equalities
\begin{gather*}
(T^+x)(t)=p_1^+(t)x(\tau_1)+(p^+-p^+_1)x(\tau_2),\quad t\in[a,b],\\
(T^-x)(t)=p_1^-(t)x(\tau_1)+(p^--p^-_1)x(\tau_2),\quad t\in[a,b],
\end{gather*}
where $p_1^+$ and $p_1^-$ are the positive and negative parts of the function $p_1$ ($p_1^+=(|p_1|+p_1)/2$, $p_1^-=(|p_1|-p_1)/2$).
Then operators $T^+$, $T^-$ satisfy equalities \eqref{e-7} and, moreover, problem  \eqref{e-3} has the same non-trivial solution as well as problem \eqref{e-8}. So, problem \eqref{e-1}--\eqref{e-2} is not uniquely solvable.
\endproof

\begin{lemma}\label{l-3}
Let non-negative numbers $\calT^+$, $\calT^-$ be given. For problem \eqref{e-8}
to have only the trivial solution for all $\tau_1$, $\tau_2\in[a,b]$ and for all functions $p_1$, $p_2\in\spaceLI$ such that
\begin{equation}\label{e-10}
    p_1(t)+p_2(t)=\calT^+-\calT^-,\quad -\calT^-\le p_1(t) \le \calT^+,\quad t\in[a,b],
\end{equation}
it is necessary and sufficient that inequalities \eqref{e-5}, \eqref{e-6} to be valid.
\end{lemma}

\begin{remark} In Lemmas \ref{l-2} and \ref{l-3}, it is sufficient to consider problem \eqref{e-8} only for the points  $\tau_1$, $\tau_2\in[a,b]$ such that $\tau_1\le\tau_2$.
\end{remark}

\proof
Let us get necessary and sufficient conditions for all problems \eqref{e-8} provided \eqref{e-9} have only the trivial solution. For any solution $y$ we have
\begin{equation*}
    y(t)=y(\tau_1)\int_a^t p_1(s)\,ds+ y(\tau_2)\int_a^t p_2(s)\,ds,\quad t\in[a,b].
\end{equation*}
Therefore, the system of equations
\begin{gather}\label{e-11}
   \left\{
   \begin{array}{l}
   C_1=C_1\int_a^{\tau_1} p_1(s)\,ds+ C_2\int_a^{\tau_1} p_2(s)\,ds,\\[5pt]
    C_2=C_1\int_a^{\tau_2} p_1(s)\,ds+ C_2\int_a^{\tau_2} p_2(s)\,ds
   \end{array}
   \right.
\end{gather}
has a solution $C_1=y(\tau_1)$, $C_2=y(\tau_2)$. Conversely, the solution
\begin{equation*}
    x(t)=C_1\int_a^t p_1(s)\,ds+ C_2\int_a^t p_2(s)\,ds,\quad t\in[a,b],
\end{equation*}
of the Cauchy problem \eqref{e-8} corresponds to
 every solution $(C_1,C_2)$ of system \eqref{e-11}. Thus, problem \eqref{e-8} has no non-trivial solutions if and only if algebraic system \eqref{e-11} has no non-trivial solutions with respect to the variables $C_1$, $C_2$,  that is, if
\begin{equation*}
    \tri\equiv\left|
                \begin{array}{cc}
                  1-\int_a^{\tau_1} p_1(s)\,ds & -\int_a^{\tau_1} p_2(s)\,ds \\[5pt]
                  -\int_a^{\tau_2} p_1(s)\,ds & 1-\int_a^{\tau_2} p_2(s)\,ds \\
                \end{array}
              \right|
              \ne 0.
\end{equation*}
Consider the determinant $\tri$ for $a\le\tau_1\le\tau_2\le b$ and for all functions $p_1$, $p_2$ satisfying conditions \eqref{e-10}. We have
\begin{equation*}
\begin{split}
\tri=    \left|
       \begin{array}{cc}
         1-\int_a^{\tau_1} p_1(s)\,ds & 1-(\a-\b)\,(\tau_1-a) \\[5pt]
         -\int_a^{\tau_2} p_1(s)\,ds & 1-(\a-\b)\,(\tau_2-a) \\
       \end{array}
     \right|
=\\[7pt]
=
   \left|
       \begin{array}{cc}
         1-\alpha      & 1-(\a-\b)\,(\tau_1-a) \\[5pt]
         -\alpha-\beta & 1-(\a-\b)\,(\tau_2-a) \\
       \end{array}
     \right|=
\\[7pt]
=
  \left|
       \begin{array}{cc}
         1-\alpha      & 1-(\a-\b)\,(\tau_1-a) \\[5pt]
         -1-\beta & -(\a-\b)\,(\tau_2-\tau_1)\\
       \end{array}
     \right|,
\end{split}
\end{equation*}
where under conditions \eqref{e-10} the values
\begin{equation*}
\alpha\equiv\int_a^{\tau_1} p_1(s)\,ds,\quad \beta\equiv\int_{\tau_1}^{\tau_2} p_1(s)\,ds,
\end{equation*}
can take arbitrary numbers from the following intervals
\begin{equation}\label{e-12}
    -(\tau_1-a)\,\b      \le \alpha \le (\tau_1-a)\,\a,\quad
    -(\tau_2-\tau_1)\,\b \le \beta  \le (\tau_2-\tau_1)\,\a.
\end{equation}

Since the determinant $\tri$ is continuous with respect to $\alpha$, $\beta$, $\tau_1$, $\tau_2$  on the connected admissible set of these parameters and $\tri=1$ for admissible values $\alpha=0$, $\beta=0$, $\tau_1=\tau_2=0$, then for all problems \eqref{e-8} to have only the trivial solution it is necessary and sufficient that $\tri$ be positive for all admissible values of these parameters. Find the minimal value of $\tri$ for fixed $\a$, $\b$ and for all rest parameters satisfying inequalities \eqref{e-12}. All problems \eqref{e-8} provided conditions \eqref{e-10} have only the trivial solutions if and only if for the pair $(\a,\b)$  this minimum is positive.

Find pairs of non-negative numbers $(\a,\b)$ such that $M\equiv\min\tri>0$, where the minimum is taken over all  $\tau_1$, $\tau_2$ such that $a\le\tau_1\le\tau_2\le b$ and over all $\alpha$, $\beta$ satisfying inequalities \eqref{e-12}. Consider the case $\b\le\a$. For $\beta=0$, $\tau_1=\tau_2$ we have $\tri=1-(\tau_1-a)\,(\a-\b)$. Then $M>0$ if and only if $(b-a)\,(\a-\b)<1$. If this inequality holds, the determinant $\tri$ takes its minimum value $M=1-(b-a)\,\a$ at $\alpha=-(\tau_1-a)\,\b$, $\beta=-(\tau_2-\tau_1)\,\b$, $\tau_2=b$, $\tau_1=a$.
Then $M>0$ if and only if  $(b-a)\,\a<1$.

In the case $\a<\b$, the minimal value $M$ is taken at
$\alpha=(\tau_1-a)\a$, $\beta=-(\tau_2-\tau_1)\b$, $\tau_2=b$,
\begin{equation*}
    \tau_1-a=\left\{
           \begin{array}{l}
           \dfrac{b-a}{2}-\dfrac{\b}{2((\calT^-)^2-(\calT^+)^2)}\ \text{ if }\
           \dfrac{\b}{(\calT^-)^2-(\calT^+)^2}<b-a; \\
           0\ \text{ if }\
           \dfrac{\b}{(\calT^-)^2-(\calT^+)^2}\ge b-a. \\
           \end{array}
           \right.
\end{equation*}
For short it is convenient  to use new variables
    $\calA\equiv (b-a)\,\a$, $\calB\equiv (b-a)\,\b$.
Then we have

\begin{equation*}
    M=\left\{
           \begin{array}{l}
           1-\calA\ \text{ if }\
           \calB\le \dfrac{1+\sqrt{1+4\calA^2}}{2}; \\[5pt]
           \dfrac{(2\calB-\calA)^2-\calA^2-
           \left(\calB^2-\calA^2-\calB+2\calA\right)^2}{4(\calB^2-\calA^2)}
           \ \text{ if }\
           \calB> \dfrac{1+\sqrt{1+4\calA^2}}{2}.
           \end{array}
           \right.
\end{equation*}
Therefore, the minimal value $M$ is positive if and only if inequalities \eqref{e-5} and  \eqref{e-6} hold.
\endproof

\begin{proof}[Proofs of Theorem \ref{t-2} and Corollary \ref{c-1}]

Assertion 1) of Corollary \ref{c-1} follows from lemmas \ref{l-1}, \ref{l-2}, and \ref{l-3}.

If in the proof of lemma \ref{l-3} not to minimize with respect to the variable  $\tau_1$, then we direct obtain the condition of the positiveness for the minimal value $M$ that is inequality \eqref{e-13} of Theorem \ref{t-2}. If we solve \eqref{e-13} with respect to the variable $\calB$, then we have condition 2) of Corollary \ref{c-1}, if with respect to the variable $\calA$, then we obtain condition 3).
\end{proof}


\begin{thebibliography}{10}
\bibitem{Hakl}
E. Bravyi, R. Hakl, A. Lomtatidze   Optimal conditions on
unique solvability of the Cauchy problem for the first order
linear functional differential equations~//
Czechoslovak Mathematical Journal. V.~52(127), No.~3.  2002. P.~513--530.

\bibitem{HL2002}
R. Hakl, A. Lomtatidze
On the Cauchy problem for first order linear
differential equations with a deviating argument~//
Archivum Mathematicum. V.~38. 2002. P.~61--71.

\bibitem{HLP2002}
R. Hakl, A. Lomtatidze, B. P\r{u}\v{z}a
New optimal conditions
for unique solvability of the Cauchy problem for first order
linear functional differential equations~// Mathematica Bohemica. V.~127, No.4. 2002. P.~509--524.

\bibitem{HLS2002-5}
R. Hakl, A. Lomtatidze, J. \v{S}remr
Some Boundary Value Problems For First
Order Scalar Functional Differential Eequations.
Folia Facult. Scien. Natur. Masar. Brunensis. Mathematica, 10. Brno:
Masaryk University, 2002. 299~p.

\bibitem{AMR} N.V. Azbelev, V.P. Maksimov, L.F. Rakhmatullina
The Elements of the Contemporary Theory of Functional Differential Equations. Methods and Applications, Institute of Computer- Assisted Studies
Moscow: Institute of Computer Researching, 2002. 384~p.

\bibitem{Zabrejko}
P.P. Zabrejko  Integral Equations. Moscow: Nauka, 1966.
448~p. (In Russian)

\bibitem{Maksimov1}
V.P. Maksimov  The Noetherian property of the  general boundary value problem for a linear functional differential equation~//
Differential equatons. V.~10, no~12. 1974. P.~2288--2291. (In Russian) 

\end{thebibliography}
\end{document}